\documentclass[12pt]{article}
\usepackage{mathrsfs}
\usepackage{graphicx}
\usepackage{lscape}
\usepackage{amsmath}
\usepackage{amssymb}
\usepackage{amsbsy}

\usepackage{amsfonts}%,srcltx}

\newtheorem{lem}{Lemma}[section]%
\newtheorem{theorem}[lem]{Theorem}%
\newtheorem{cor}[lem]{Corollary}%
\newtheorem{prop}[lem]{Proposition}%

\def\a{\alpha} \def\b{\beta} \def\g{\gamma} \def\d{\delta}

\def\G{\Gamma}

 \def\lg{\langle} \def\rg{\rangle}
\def\nd{\mathrel{\bigm|\kern-.7em/}}

\def\f{\noindent}

\def\Aut{\hbox{\rm Aut}}

\def\Cay{\hbox{\rm Cay}}

\def\Sym{\hbox{\rm Sym}}
\def\Alt{\hbox{\rm Alt}}

\def\GP{\hbox{\rm GP}}
\def\Q{\hbox{\rm Q}}

\def\demo{\f {\bf Proof.}\hskip10pt}

\newcommand{\qed}{\mbox{\raisebox{0.7ex}{\fbox{}}} \vspace{4truemm}}
\def\mz{{\mathbb Z}}

\def\K{{\bf K}}
\def\1{{\bf 1}}

\begin{document}

\title{Every finite group has a normal bi-Cayley graph}
\vspace{4 true mm}

%\footnotetext{\footnotesize{\em E-mail}: chenghw2002$@$sina.com,\
%jxzhou$@$bjtu.edu.cn}\vspace{0.5mm}

\author {{\sc Jin-Xin Zhou}\\[1ex]
{\small\em Mathematics, Beijing Jiaotong University,
Beijing 100044, P.R. China}\\
{\small{\em E-mail}: \texttt{jxzhou@bjtu.edu.cn}} }
\date{}
\maketitle

\begin{abstract}
A graph $\G$ with a group $H$ of automorphisms acting semiregularly on the vertices with two orbits is called a {\em bi-Cayley graph} over $H$. When  $H$ is a normal subgroup of $\Aut(\G)$, we say that $\G$ is {\em normal} with respect to $H$. In this paper, we show that every finite group has a connected normal bi-Cayley graph. This improves Theorem~5 of [M. Arezoomand, B. Taeri, Normality of 2-Cayley digraphs, Discrete Math. 338 (2015) 41--47], and provides a positive answer to the Question of the above paper.

\bigskip

\end{abstract}
\thispagestyle{empty}

\section{Introduction}
Throughout this paper, groups are assumed to be finite, and graphs are assumed to be finite, connected, simple and undirected.
For the group-theoretic and graph-theoretic terminology not defined here we refer the
reader to \cite{BMBook,WI}.

Let $G$ be a permutation group on a set $\Omega$ and $\a\in \Omega$. Denote by $G_\a$ the stabilizer of $\a$ in $G$, that is, the subgroup of $G$ fixing the point $\a$. We say that $G$ is {\em semiregular} on $\Omega$ if $G_\a=1$ for every $\a\in \Omega$ and {\em regular} if $G$ is transitive and semiregular. It is well-known that a graph $\G$ is a {\em Cayley graph} if it has an automorphism group acting regularly on its vertex set (see \cite[Lemma~16.3]{B}). If we, instead, require the graph $\G$ having a group of automorphisms acting semiregularly on its vertex set with two orbits, then we obtain a so-called {\em bi-Cayley graph}.

Cayley graph is usually defined in the following way. Given a finite group $G$ and an inverse closed subset $S\subseteq G\setminus\{1\}$, the {\em Cayley graph} $\Cay(G,S)$ on $G$ with respect to $S$ is a graph with vertex set $G$ and edge set $\{\{g,sg\}\mid g\in G,s\in S\}$. For any $g\in G$, $R(g)$ is the permutation of $G$ defined by $R(g): x\mapsto xg$ for $x\in G$. %The {\em right regular representation} of $G$ is the subgroup of Sym$(G)$ defined by
Set $R(G):=\{R(g)\ |\ g\in G\}$. It is well-known that $R(G)$ is a subgroup of $\Aut(\Cay(G,S))$. In 1981, Godsil~\cite{Godsil1981} proved that the normalizer of $R(G)$ in $\Aut(\Cay(G,S))$ is $R(G)\rtimes\Aut(G,S)$, where $\Aut(G,S)$ is the group of automorphisms of $G$ fixing the set $S$ set-wise. This result has been successfully used in characterizing various families of Cayley graphs $\Cay(G,S)$ such that $R(G)=\Aut(\Cay(G,S))$ (see, for example, \cite{Godsil1981,Godsil1983}). A {Cayley graph} $\Cay(G,S)$ is said to be {\em normal} if $R(G)$ is normal in $\Aut(\Cay(G,S))$. This concept was introduced by Xu in \cite{X1}, and for more results about normal Cayley graphs, we refer the reader to \cite{FengLuXu}.

Similarly, every bi-Cayley graph admits the following concrete realization. Given a group $H$, let $R$, $L$ and $S$ be subsets of $H$ such that $R^{-1}=R$, $L^{-1}=L$ and $R\cup L$ does not contain the identity element of $H$. The {\em bi-Cayley graph} over $H$ relative to the triple $(R, L, S)$, denoted by BiCay($H$,~$R$,~$L$,~$S$), is the graph having vertex set the union of the right part $H_{0}=\{h_{0}~|~h\in H\}$ and the left part $H_{1}=\{h_{1}~|~h\in H\}$, and edge set the union of the right edges $\{\{h_{0},~g_{0}\}~|~gh^{-1}\in R\}$, the left edges $\{\{h_{1},~g_{1}\}~|~gh^{-1}\in L\}$ and the spokes $\{\{h_{0},~g_{1}\}~|~gh^{-1}\in S\}$.
Let $\G=$BiCay$(H, R, L, S)$. For $g\in H$, define a permutation $BR(g)$ on the vertices of $\G$ by the rule
$$h_i^{BR(g)}=(hg)_i, \forall i\in\mz_2, h\in H.$$
Then $BR(H)=\{BR(g)\ |\ g\in H\}$ is a semiregular subgroup of $\Aut(\G)$ which is isomorphic to $H$ and has $H_0$ and $H_1$ as its two orbits.
When $BR(H)$ is normal in $\Aut(\G)$, the bi-Cayley graph $\G=$BiCay$(H, R, L, S)$ will be called a {\em normal bi-Cayley graph} over $H$ (see \cite{AT-normal-bi-Cay} or \cite{ZF-auto}).

Wang et al. in \cite{WWX} determined the groups having a connected normal Cayley graph.

\begin{prop}\label{prop-normal-cay}
Every finite group $G$ has a normal Cayley graph unless $G\cong C_4\times C_2$ or $G\cong\Q_8\times C_2^r(r\geq0)$.
\end{prop}

Following up this result, Arezoomand and Taeri in \cite{AT-normal-bi-Cay} asked: Which finite groups have normal bi-Cayley graphs? They also gave a partial answer to this question by proving that every finite group $G\not\cong\Q_8\times C_2^r(r\geq0)$ has at least one normal bi-Cayley graph. At the end of \cite{AT-normal-bi-Cay}, the authors asked the following question:\medskip

\f{\bf Question~A}\ {\rm\cite[Question]{AT-normal-bi-Cay}}\
Is there any normal bi-Cayley graph over $G\not\cong\Q_8\times C_2^r$ for each $r\geq0$?\medskip

%\begin{prop}{\rm\cite[Theorem~5]{AT-normal-bi-Cay}}\label{prop-normal-bi-cay}
%Every finite group $G\not\cong\Q_8\times\mz_2^r(r\geq0)$ has at least one normal bi-Cayley graph.
%\end{prop}

We remark that for every finite group $G\not\cong\Q_8\times C_2^r(r\geq0)$, the normal bi-Cayley graph over $G$ constructed in the proof of \cite[Theorem~5]{AT-normal-bi-Cay} is not of regular valency, and so is not vertex-transitive. So it is natural to ask the following question.\medskip

\f{\bf Question~B}\
Is there any vertex-transitive normal bi-Cayley graph over a finite group $G$?\medskip

In this paper, Questions~A and B are answered in positive. The following is the main result of this paper.

\begin{theorem}\label{th-main}
Every finite group has a vertex-transitive normal bi-Cayley graph.
\end{theorem}

To end this section we give some {notation} which will be used in this paper.
For a positive integer $n$, denote by $C_n$ the cyclic group of order $n$, by $\mz_n$ the ring of integers modulo $n$,  by $D_{2n}$ the dihedral group of order $2n$, and by $\Alt(n), \Sym(n)$ the alternating group and symmetric group of degree $n$, respectively. Denote by $\Q_8$ the quaternion group. For two groups $M$ and $N$, $N\rtimes M$ denotes a semidirect product of $N$ by $M$. Let $G$ be a finite group. Denote by $1_G$ the identity element of $G$. 

For a finite, simple and undirected graph $\G$, we use $V(\G)$, $E(\G)$ and $\Aut(\G)$ to denote its vertex set, edge set and full automorphism group, respectively, and for any $u, v\in V(\G)$, $u\sim v$ means $u$ and $v$ are adjacent. A graph $\G$ is said to be {\em vertex-transitive} if its full automorphism group $\Aut(\G)$ acts transitively on its vertex set For any subset $B$ of $V(\G)$, the subgraph of $\G$ induced by $B$
will be denoted by $\G[B]$. %, and the neighborhood of $B$ in $\G$ is defined as $N_\G(B)=\bigcup_{v\in B}(\G(v))-B$.

\section{Cartesian products}

The {\em Cartesian product} $X\square Y$ of graphs $X$ and $Y$ is a graph with vertex set $V(X)\times V(Y)$, and vertices $(u, x)$ and $(v, y)$ adjacent if and only if $u=v$ and $x\sim y$ in $Y$, or $x=y$ and $u\sim v$ in $X$.

A non-trivial graph $X$ is {\em prime} if it is not isomorphic to a Cartesian product of two smaller graphs. The following proposition shows the uniqueness of the prime factor decomposition of connected graphs with respect to the Cartesian product.

\begin{prop}{\rm\cite[Theorem~6.6]{Imrich-Klavzar}}\label{decomposition}
Every connected finite graph can be decomposed as a Cartesian product of prime graphs, uniquely up to isomorphism and the order of the factors.
\end{prop}

Two non-trivial graphs are {\em relatively prime} (w.r.t. Cartesian product) if they have no non-trivial common factor. Now we consider the automorphisms of Cartesian product of graphs.

\begin{prop}{\rm\cite[Theorem~6.10]{Imrich-Klavzar}}\label{auto-decomposition}
Suppose $\phi$ is an automorphism of a connected graph $\G$ with prime factor decomposition $\G=\G_1\square \G_2\square \cdots\square \G_k$. Then there is a permutation $\pi$ of $\{1,2,\ldots,k\}$ and isomorphisms $\phi_i: \G_{\pi(i)}\rightarrow \G_i$ for which
$$\phi(x_1, x_2, \ldots, x_k)=(\phi_1(x_{\pi(1)}), \phi_2(x_{\pi(2)}), \ldots, \phi_k(x_{\pi(k)})).$$
\end{prop}

Let $\G$ be a connected graph  with prime factor decomposition $\G=\G_1\square \G_2\square \cdots\square \G_k$. Take $(x_1, x_2, \ldots, x_k)\in V(\G)$. For each $1\leq i\leq k$, set $$V_{x_{i}}=\{(x_1, \ldots, x_{i-1}, y, x_{i+1}, \ldots, x_k)\ |\ y\in V(\G_i)\}.$$
For any $\phi\in\Aut(\G)$, suppose $V_{x_{i}}^\phi\cap V_{x_i}\neq\emptyset$. Take $(x_1, \ldots, x_{i-1}, y, x_{i+1}, \ldots, x_k)\in V_{x_{i}}^\phi\cap V_{x_i}$. Then there exists $y'\in V(\G_i)$ such that
\begin{equation}\label{eq-1}
\phi(x_1, \ldots, x_{i-1}, y', x_{i+1}, \ldots, x_k)=(x_1, \ldots, x_{i-1}, y, x_{i+1}, \ldots, x_k).
\end{equation}
By Proposition~\ref{auto-decomposition}, there is a permutation $\pi$ of $\{1,2,\ldots,k\}$ and isomorphisms $\phi_i: \G_{\pi(i)}\rightarrow \G_i$ for which
\begin{equation}\label{eq-2}
\phi(x_1, \ldots, x_{i-1}, x_i=y', x_{i+1}, \ldots, x_k)=(\phi_1(x_{\pi(1)}), \phi_2(x_{\pi(2)}), \ldots, \phi_k(x_{\pi(k)})).
\end{equation}
Combining Equation~\ref{eq-1} with Equation~\ref{eq-2}, we have the permutation $\pi$ is the identity, and each $\phi_j$ is an automorphism of $\G_j$, and if $j\neq i$, then $\phi_j(x_j)=x_j$. Consequently, we have $V_{x_i}^{\phi}=V_{x_i}$, and so $V_{x_i}$ is a block of imprimitivity of $\Aut(\G)$.

\begin{cor}\label{prop-cartersian}
Use the same notation given in the above paragraph. Then the following hold.
\begin{enumerate}
  \item [{\rm (1)}]\ $V_{x_i}$ is a block of imprimitivity of $\Aut(\G)$.
  \item [{\rm (2)}]{\rm(\cite[Corollary~6.12]{Imrich-Klavzar})}\ If $\G_1, \G_2, \ldots, \G_k$ are relatively prime, then $\Aut(\G)=\Aut(\G_1)\times\Aut(\G_2)\times\cdots\times\Aut(\G_k)$.
\end{enumerate}
\end{cor}

The following theorem provides a method of constructing normal bi-Cayley graphs.

\begin{theorem}\label{lem-cartersian}
Let $X$ be a connected normal bi-Cayley graph over a group $H$, and let $Y$ be a connected normal Cayley graph over a group $K$. If $X$ and $Y$ are relatively prime, then $X\square Y$ is also a normal bi-Cayley graph over the group $H\times K$.
\end{theorem}

\demo Assume that $X$ and $Y$ are relatively prime. By Corollary~\ref{prop-cartersian}, $\Aut(X\square Y)=\Aut(X)\times\Aut(Y)$. Since $X$ is a connected normal bi-Cayley graph over $H$, one has $BR(H)\unlhd\Aut(X)$, and since $Y$ is a connected normal Cayley graph over a group $K$, one has $R(K)\unlhd\Aut(Y)$. Then $BR(H)\times R(K)$ is a normal subgroup of $\Aut(X\square Y)=\Aut(X)\times\Aut(Y)$. Note that $BR(H)$ acts semiregularly on $V(X)$ with two orbits, and $R(K)$ acts regularly on $V(Y)$. It follows that $BR(H)\times R(K)$ acts semiregularly on $V(X)\times V(Y)$ with two orbits, and thereby $X\square Y$ is also a normal bi-Cayley graph over the group $H\times K$. \hfill\qed

\section{Normal bi-Cayley graphs over ${\rm Q_8}\times C_2^{r} (r\geq 0)$}

In this section, we shall answer Question~A in positive.

\subsection{The $n$-dimensional hypercube}

For $n\geq 1$, the {\em $n$-dimensional hypercube}, denoted by $Q_n$, is the graph whose vertices are all the $n$-tuples of $0$'s and $1$'s with two $n$-tuples being adjacent if and only if they differ in exactly one place.

Let $N= C_2^n$ be an elementary abelian $2$-group of order $2^n$ with a minimum generating set $S=\{s_1, s_2, s_3, \ldots, s_n\}$. By the definition of $Q_n$, we have $\Cay(N, S)\cong Q_n$. For convenience of the statement, we assume that $Q_n=\Cay(N, S)$. If $n=1$, then $Q_1=\K_2$ and so $\Aut(Q_1)=N$. In what follows, assume that $n\geq 2$.  It is easy to observe that for any distinct $s_i, s_j$ there is a unique $4$-cycle in $Q_n$ passing through $1_N, s_i, s_j$, where $1_N$ is the identity element of $N$. So if a subgroup of $\Aut(Q_n)$ fixes $S$ pointwise, then it also fixes every vertex at distance $2$ from $1_N$. By the connectedness and vertex-transitivity of $Q_n$, we have $\Aut(Q_n)_{1_N}$ acts faithfully on $S$. It follows that $\Aut(Q_n)_{1_N}\lesssim\Sym(n)$. On the other hand, it is easy to see that each permutation on $S$  induce an automorphism of $N$, and so $\Aut(N, S)\cong\Sym(n)$. Since $\Aut(N, S)\leq\Aut(Q_n)_{1_N}$, one has $\Aut(Q_n)_{1_N}=\Aut(N, S)\cong\Sym(n)$. Consequently, we have $\Aut(Q_n)=R(N)\rtimes\Aut(N, S)\cong N\rtimes\Sym(n)$.

Note that $Q_n$ is bipartite. Let $\Aut(Q_n)^*$ be the kernel of $\Aut(Q_n)$ acting on the two partition sets of $Q_n$. Let $E=R(N)\cap \Aut(Q_n)^*$. Then $E\unlhd \Aut(Q_n)^*$ and $E\unlhd R(N)$. It follows that $E\unlhd \Aut(Q_n)^*R(N)=\Aut(Q_n)$. Clearly, $E$ acts semiregularly on $V(Q_n)$ with two orbits. Thus, we have the following lemma.

\begin{lem}\label{lem-q}
Use the same notation as in the above three paragraphs. For any $n\geq 1,$ $Q_n$ is a normal Cayley graph over $N$, and $Q_n$ is also a normal bi-Cayley graph over $E$.
\end{lem}

\subsection{The M\"obius-Kantor graph}

\begin{figure}[ht]
\begin{center}
\unitlength 4mm
\begin{picture}(15,12)

\put(5,1){\circle{0.35}} \put(4.5,0){$1$} \put(5,1){\line(1,0){4}}
\put(5,1){\vector(-1,1){1.5}} \put(5,1){\vector(1,0){2}}
\put(7,0.4){$z_2$} \put(5,1){\line(3,5){1}}
\put(5,1){\line(-1,1){2.5}}

\put(9,1){\circle*{0.35}}\put(9.5,0){$2$} \put(9,1){\line(-3,5){1}}
\put(9,1){\line(1,1){2.5}}\put(9,1){\vector(1,1){1.5}}
\put(10.2,1.7){$z_3$}

\put(2.5,3.5){\circle*{0.35}}\put(1.8,3.5){$0$}
\put(2.5,3.5){\line(0,1){4}}\put(2.5,3.5){\line(3,1){2}}
%\put(2.5,3.5){\vector(1,-1){1.5}}
\put(3,1.8){$z_9$}

\put(11.5,3.5){\circle{0.35}}\put(12,3.5){$3$}
\put(11.5,3.5){\line(0,1){4}}\put(11.5,3.5){\line(-3,1){2}}
\put(11.5,3.5){\vector(0,1){2}} \put(11.6,5.5){$z_4$}

\put(2.5,7.5){\circle{0.35}} \put(1.8,7.5){$7$}
\put(2.5,7.5){\line(1,1){2.5}}\put(2.5,7.5){\line(3,-1){2}}
\put(2.5,7.5){\vector(0,-1){2}} \put(1.6,5.5){$z_8$}

\put(11.5,7.5){\circle*{0.35}} \put(12,7.5){$4$}
\put(11.5,7.5){\line(-1,1){2.5}}\put(11.5,7.5){\line(-3,-1){2}}
\put(11.5,7.5){\vector(-1,1){1.5}} \put(10,9){$z_5$}

\put(5,10){\circle*{0.35}} \put(5,10.5){$6$}
\put(5,10){\line(1,0){4}}\put(5,10){\line(3,-5){1}}
\put(5,10){\vector(-1,-1){1.5}} \put(3.2,9){$z_7$}

\put(9,10){\circle{0.35}} \put(9,10.5){$5$}
\put(9,10){\line(-3,-5){1}} \put(9,10){\vector(-1,0){2}}
\put(7,10.2){$z_6$}

\put(6,2.6){\circle*{0.35}}\put(6.2,2){$1'$}
\put(6,2.6){\line(0,1){5.7}}\put(6,2.6){\line(5,6){3.5}}

\put(8,2.6){\circle{0.35}}\put(8.3,2.3){$2'$}
\put(8,2.6){\line(0,1){5.7}}\put(8,2.6){\line(-5,6){3.5}}

\put(6,8.3){\circle{0.35}}\put(6.3,8.3){$6'$}
\put(6,8.3){\line(5,-6){3.5}}

\put(8,8.3){\circle*{0.35}}\put(8.35,8.2){$5'$}
\put(8,8.3){\line(-5,-6){3.5}}

\put(9.5,4.2){\circle*{0.35}}\put(9.8,4.2){$3'$}
\put(9.5,4.2){\line(-1,0){5}} \put(9.5,4.2){\vector(-1,0){2.7}}
\put(6.7,4.4){$z_1$}

\put(9.5,6.8){\circle{0.35}}\put(9.8,6.2){$4'$}
\put(9.5,6.8){\line(-1,0){5}}

\put(4.5,6.8){\circle*{0.35}}\put(4.2,7.1){$7'$}

\put(4.5,4.2){\circle{0.35}}\put(4.2,4.5){$0'$}
\end{picture}
\end{center}\vspace{-.5cm}
\caption{The M\"obius-Kantor graph $\GP(8, 3)$} \label{gp}
\end{figure}
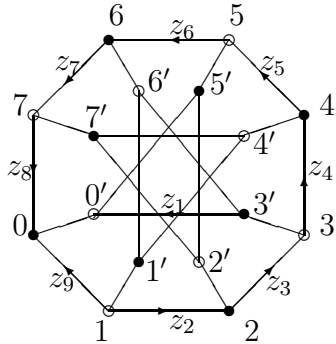

The M\"obius-Kantor graph $\GP(8, 3)$ is a graph with
vertex set $V=\{i,i'\ |\ i\in\mz_8\}$ and edge set the union of the
{\em outer edges} $\{\{i,i+1\}\ |\ i\in \mz_8 \}$, the {\em inner
edges} $\{\{i',(i+3)'\}\ |\ i\in\mz_8\}$, and the {\em spokes}
$\{\{i,i'\}\ |\ i\in \mz_8\}$ (see Fig.~\ref{gp}). Note that $\GP(8, 3)$ is a bipartite
graph with bipartition sets $B_1=\{1,3,5,7,0',2',4',6'\}$ and
$B_2=\{0,2,4,6,1',3',5',7'\}$.

In \cite{ZF2012}, the edge-transitive automorphism groups of $\Aut(\GP(8,3))$ were determined. We first introduce the following automorphisms of $\GP(8,3)$, represented as permutations on
the vertex set $V$:

$$%
\begin{array}{lll}
 \a&=&(1\ 3\ 5\ 7)(0\ 2\ 4\ 6)(1'\ 3'\ 5'\ 7')(0'\ 2'\ 4'\ 6'), \\
   \b&=&(0\ 1'\ 2)(0'\ 6'\ 3)(4\ 5'\ 6)(7\ 4'\ 2'), \\
   \g&=&(1\ 1')(2\ 6')(3\ 3')(4\ 0')(5\ 5')(6\ 2')(7\ 7')(0\ 4'),\\
\d&=&(1\ 1')(2\ 4')(3\ 7')(4\ 2')(5\ 5')(6\ 0')(7\ 3')(0\ 6').
\end{array}
$$%

By \cite[Lemma~3.1]{ZF2012}, we have $\lg\a,\b\rg=\lg\a,\a^\b\rg\rtimes\lg\b\rg\cong{\rm Q}_8\rtimes\mz_3$, where ${\rm Q}_8$ is the quaternion group, and moreover, $\lg\a,\b\rg\unlhd\Aut(\GP(8,3))$. Clearly, $\lg\a,\a^\b\rg\cong{\rm Q}_8$ is the Sylow $2$-subgroup of $\lg\a,\b\rg$, so $\lg\a,\a^\b\rg$ is characteristic in $\lg\a,\b\rg$, and then it is normal in $\Aut(\GP(8,3))$ because $\lg\a,\b\rg\unlhd\Aut(\GP(8,3))$. For convenience of the statement, we let ${\rm Q}_8=\lg \a,\a^\b\rg$. It is easy to see that ${\rm Q}_8$ acts semiregularly on $V$ with two orbits $B_1$ and $B_2$. Thus we have the following lemma.

\begin{lem}\label{lem-gp}
$\GP(8,3)$ is a normal bi-Cayley graph over ${\rm Q}_8$.
\end{lem}

\subsection{An answer to Question~A}

Noting that $\GP(8, 3)$ is of girth $6$, $\GP(8, 3)$ is prime. For each $r\geq 1$, it is easy to see that $Q_{r}=\underbrace{\K_2\square\K_2\square\cdots\square\K_2}_{n\ {\rm times}}$. So, $Q_n$ and $\GP(8, 3)$ are relatively prime. Now combining Lemmas~\ref{lem-q} and \ref{lem-gp} with Theorem~\ref{lem-cartersian}, we can obtain the following theorem.

\begin{theorem}\label{th-answer}
For each $r\geq 1$, $\GP(8, 3)\times Q_r$ is a vertex-transitive normal bi-Cayley graph over ${\rm Q}_8\times N$, where $E\cong C_2^{r}$.
\end{theorem}

\section{Proof of Theorem~\ref{th-main}}

The proof of Theorem~\ref{th-main} will be completed by the following lemmas. Let $G$ be a group. A Cayley graph $\G=\Cay(G, S)$ on $G$ is said to be a {\em GRR} of $G$ if $\Aut(\G)=R(G)$.

\begin{lem}\label{th-general-GRR}
Let $G$ be a group admitting a GRR $\G$. Then $\G\square \K_2$ is a normal bi-Cayley graph over the group $G$.
\end{lem}

\demo If $\K_2$ and $\G$ are relatively prime, then by Corollary~\ref{prop-cartersian}~(2), we have $\Aut(\G\square\K_2)=\Aut(\G)\times\Aut(\K_2)$.
Clearly, $R(G)\times {\bf 1}$ acts semiregularly on $V(\G\square\K_2)$ with two orbits, and $R(G)\times {\bf 1}\unlhd \Aut(\G\square\K_2)$, where ${\bf 1}$ is the identity of $\Aut(\K_2)$. It follows that $\G\square\K_2$ is a normal bi-Cayley graph over the group $G$.

Suppose that $\K_2$ is also a prime factor of $\G$. Let $\G=\G_1\square\underbrace{\K_2\square\cdots\square\K_2}_{m\ {\rm times}}$ be such that $\G_1$ is coprime to $\K_2$. From Corollary~\ref{prop-cartersian}~(2) it follows that $G=\Aut(\G)=\Aut(\G_1)\times\Aut(\K_2\square\cdots\square\K_2)$. Since $\G$ is a GRR of $G$, one has $m=1$, and so $\G\square\K_2=\G_1\square\K_2\square\K_2$. Then $G=\Aut(\G_1)\times \Aut(\K_2)$, and $\G_1$ is a GRR of $\Aut(\G_1)$. By Lemma~\ref{lem-q}, $\K_2\square\K_2$ is a normal bi-Cayley graph over $C_2$, and by Theorem~\ref{lem-cartersian}, $\G\square\K_2$ is a normal bi-Cayley graph over $\Aut(\G_1)\times C_2\cong G$.\hfill\qed

A group $G$ is called {\em generalized dicyclic group} if it is non-abelian and has an abelian subgroup $L$ of index $2$ and an element $b\in G\setminus L$ of order $4$ such that $b^{-1}xb=x^{-1}$ for every $x\in L$.

The following theorem gives a list of groups having no GRR (see \cite{Godsil1978}).

\begin{theorem}\label{GRR}
Every finite group $G$ admits a GRR unless $G$ belongs to one of the following classes of groups:
\begin{enumerate}
  \item [{\rm (I)}]\ Class C: abelian groups of exponent greater than two;
  \item [{\rm (II)}]\ Class D: the generalized dicyclic groups;
  \item [{\rm (III)}]\ Class E: the following thirteen ``exceptional groups":
  \begin{enumerate}
    \item [{\rm (1)}]\ $\mz_2^2, \mz_2^3, \mz_2^4$;
    \item [{\rm (2)}]\ $D_6, D_8, D_{10}$;
    \item [{\rm (3)}]\ $A_4$;
    \item [{\rm (4)}]\ $\lg a,b,c \ |\ a^2=b^2=c^2=1, abc=bca=cab\rg$;
    \item [{\rm (5)}]\ $\lg a,b \ |\ a^8=b^2=1, bab=a^5\rg$;
    \item [{\rm (6)}]\ $\lg a,b,c\ |\ a^3=b^3=c^2=1, ac=ca, (ab)^2=(cb)^2=1\rg$;
    \item [{\rm (7)}]\ $\lg a,b,c\ |\ a^3=b^3=c^3=1, ac=ca, bc=cb, c=a^{-1}b^{-1}cb\rg$;
    \item [{\rm (8)}]\ ${\rm Q}_8\times \mz_3, {\rm Q}_8\times \mz_4$.
  \end{enumerate}
\end{enumerate}
\end{theorem}

\begin{lem}\label{th-general-C}
Let $G$ be a group in Class~C of Theorem~\ref{GRR}. Then $G$ has a normal bi-Cayley graph.
\end{lem}

\demo Since $G$ is abelian, $G$ has an automorphism $\a$ such that $\a$ maps every element of $G$ to its inverse. Set $H=G\rtimes\lg\a\rg$.
If $H$ has a GRR $\G$, then $\G$ is also a normal bi-Cayley graph over $G$. Suppose that $H$ has no GRR. Then by Theorem~\ref{GRR} we have $H$ is one of the groups in Class E (2) and (6). By Lemma~\ref{th-general-E}, $G$ has a normal bi-Cayley graph\hfill\qed

\begin{lem}\label{th-general-D}
Let $G$ be a group in Class~D  of Theorem~\ref{GRR}. Then $G$ has a normal bi-Cayley graph.
\end{lem}

\demo If $G\cong {\rm Q_8}\times\mz_2^r$ for some $r\geq0$, then by Theorem~\ref{th-answer} and Lemma~\ref{lem-gp}, $G$ has a normal bi-Cayley graph. In what follows, we assume that $G\not\cong{\rm Q_8}\times\mz_2^r$ for any $r\geq0$.

By the definition of generalized dicyclic group, $G$ has an abelian subgroup $L$ of index $2$ and an element $b$ of order $4$ such that $G=L\lg b\rg$ and $b^{-1}ab=a^{-1}$ for any $a\in L$. We first prove the following claim.\medskip

\f{\bf Claim}\ Let $\G=\Cay(G, S)$ be a connected Cayley graph of $G$. If $\K_2$ is a factor of $\G$, then there exists an involution $s\in S$ such that $G=\lg S-s\rg\times\lg s\rg$.

Assume that $\K_2$ is a factor of $\G$. We may let $\G\cong\G'=\G_1\square\K_2$. Assume that $V(\K_2)=\{0, 1\}$. For any $v\in V(\G_1)$, by Corollary~\ref{lem-cartersian}, $V_{v}=\{(v, 0), (v, 1)\}$ is a block of imprimitivity of $\Aut(\G_1\square\K_2)$. Clearly, $V_{v}$ is also an edge of $\G'$. Let $F'=\{\{(u, 0), (u, 1)\}\ |\ u\in V(\G_1)\}$. Then $F'$ is a $1$-factor of $\G'$, and $\G-F'$ is just a union of two copies of $\G_1$.

As $\G\cong\G'$, $\G$ has an edge which is a block of imprimitivity of $\Aut(\G)$. Since $\G$ is vertex-transitive, there exists $s\in S$ such that $\{1, s\}$ is a block of imprimitivity of $\Aut(\G)$. Since $R(G)$ acts regularly on $V(\G)=G$ by right multiplication, we have $R(G)_{\{1, s\}}=\{1, s\}\cong\mz_2$. Thus $s$ is an involution. If $s\notin L$, then since $|G: L|=2$, one has $G=L\rtimes\lg s\rg$. It follows that $s=xb$ or $xb^{-1}$ for some $x\in L$. So $s^{-1}as=a^{-1}$ for each $a\in L$. This implies that every element of $G$ outside $L$ would be an involution, a contradiction. Thus, $s\in L$ and so $s$ is in the center of $G$.

Recall that $F'$ is a $1$-factor of $\G'$. Correspondingly, $F=\{\{g, sg\}\ |\ g\in G\}$ is a $1$-factor of $\G$, and $\G-F$ is disconnected with two isomorphic components. This implies that $\lg S-s\rg$ has index $2$ in $G$. Thus, $G=\lg S-s\rg\times\lg s\rg$, completing the proof of Claim.\medskip

Now let $G=G_1\times G_2$ be such that $G_2\cong C_2^r$ with $r\geq 0$ and $G_1$ can not be decomposed into a direct product of two proper subgroups one of which is an elementary abelian $2$-group. Then $|G_1: L\cap G_1|=2$ and $G_2\leq L$. Let $b=g_1g_2$ with $g_1\in G_1$ and $g_2\in G_2$. Note that $b$ has order $4$ and $g_2$ has order at most $2$. Clearly, $g_1$ commutes with $g_2$, so $g_1$ has order $4$. Furthermore, for any $a\in G_1\cap L$, we have $$a^{-1}=b^{-1}ab=(g_1g_2)^{-1}a(g_1g_2)=g_1^{-1}ag_1.$$
It follows that $G_1$ is also a generalized dicyclic group.

By Proposition~\ref{prop-normal-cay}, $G_1$ has a connected normal Cayley graph, say $\Sigma$. By Claim, $\Sigma$ is coprime to $\K_2$. Let $\Lambda=\underbrace{\K_2\square\cdots\square\K_2}_{r+1\ {\rm times}}$. By Lemma~\ref{lem-q}, $\Lambda$ is a connected normal bi-Cayley graph over $C_2^r$. By Theorem~\ref{lem-cartersian}, $\Sigma\square\Lambda$ is a connected normal bi-Cayley graph over $G_1\times C_2^r\cong G$.\hfill\qed

\begin{lem}\label{th-general-E}
Let $G$ be a group in Class~E of Theorem~\ref{GRR}. Then $G$ has a normal bi-Cayley graph.
\end{lem}

\demo By Lemma~\ref{lem-q}, each of the groups in Class~E(1) has a connected normal bi-Cayley graph.

Let $G=D_{2n}=\lg a,b\ |\ a^n=b^2=1, b^{-1}ab=a^{-1}\rg$ with $n\geq 3$. Let $\G=\Cay(G, \{ab, b\})$. Then $\G$ is a cycle of length $2n$, and so $\G$ is coprime to $\K_2$. By Theorem~\ref{lem-cartersian}, $\G\square\K_2$ is a connected normal bi-Cayley graph over $G$. Thus, each of the groups in Class~E(2) has a connected normal bi-Cayley graph.

Let $G={\rm Alt}(4)$ and let $\G=\Cay(G, \{(1\ 2\ 3), (1\ 3\ 2), (1\ 2\ 4), (1\ 4\ 2)\})$. By Magma~\cite{BCP}, we have $\G\square\K_2$ is a connected normal bi-Cayley graph over ${\rm Alt}(4)$.

Let $G=\lg a,b,c \ |\ a^2=b^2=c^2=1, abc=bca=cab\rg$ be the group in  Class~E(4). Let $\G=\Cay(G, \{a,b,c\})$. By Magma~\cite{BCP}, $\G$ is a connected trivalent normal Cayley graph over $G$ and $\G$ has girth $6$. Hence, $\G$ is coprime to $\K_2$. By Theorem~\ref{lem-cartersian}, $\G\square\K_2$ is a connected normal bi-Cayley graph over $G$.

Let $G=\lg a,b \ |\ a^8=b^2=1, bab=a^5\rg$ be the group in  Class~E(5). Let $\G=\Cay(G, \{a, a^{-1}, b, a^4, a^4b\})$. By \cite[Lemma~6]{WWX}, $\G$ is a connected normal Cayley graph over $G$, and by Magma, $\Aut(\G\square\K_2)=\Aut(\G)\times\mz_2$. Thus, $\G\square\K_2$ is a normal bi-Cayley graph over $G$.

Let $G=\lg a,b,c\ |\ a^3=b^3=c^2=1, ac=ca, (ab)^2=(cb)^2=1\rg$  be the group in  Class~E(6). Let $\G=\Cay(G, \{c, ca, cb\})$. By Magma~\cite{BCP}, $\G$ is a connected trivalent normal Cayley graph over $G$ and $\G$ has girth $6$. Hence, $\G$ is coprime to $\K_2$. By Lemma~\ref{lem-cartersian}, $\G\square\K_2$ is a connected normal bi-Cayley graph over $G$.

Let $G=\lg a,b,c\ |\ a^3=b^3=c^3=1, ac=ca, bc=cb, c=a^{-1}b^{-1}cb\rg$   be the group in  Class~E(7). Let $\G=\Cay(G, \{a,b,a^{-1},b^{-1}\})$. By Magma~\cite{BCP}, $\G$ is a connected trivalent normal Cayley graph over $G$. Since $G$ has order $27$, $\G$ is coprime to $\K_2$. By Theorem~\ref{lem-cartersian}, $\G\square\K_2$ is a connected normal bi-Cayley graph over $G$.

Finally, we consider the groups in  Class~E(8). By Lemma~\ref{lem-gp}, $\GP(8,3)$ is a normal bi-Cayley graph over ${\rm Q}_8$. For $n\geq 3$, let $\mz_n=\lg a\rg$ and let ${\rm C}_n=\Cay(\mz_n, \{a, a^{-1}\})$. Clearly, ${\rm C}_n$ is a normal Cayley graph over $\mz_n$. Since $\GP(8,3)$ is of girth $6$, $\GP(8,3)$ is coprime to ${\rm C_n}$. By Theorem~\ref{lem-cartersian}, $\GP(8,3)\square{\rm C}_n$ is a connected normal bi-Cayley graph over ${\rm Q}_8\times\mz_n$. Thus each of the groups in  Class~E~(8) has  a connected normal bi-Cayley graph.\hfill\qed

\f{\bf Proof of Theorem~\ref{th-main}}\ Let $G$ be a finite group. If $G$ has a GRR, then by Lemma~\ref{th-general-GRR}, $G$ has a connected normal bi-Cayley graph. If $G$ does not has a GRR, then the theorem follows from Lemmas~\ref{th-general-C}-\ref{th-general-E} and \ref{lem-gp}.\hfill\qed

\medskip
\f {\bf Acknowledgements:}\ This work was partially supported by the National
Natural Science Foundation of China (11271012) and the Fundamental
Research Funds for the Central Universities (2015JBM110).

{}
\end{document}